\documentclass[12pt]{amsart}
\usepackage[cp1251]{inputenc}
\usepackage[english]{babel}
\usepackage{amsmath}
\usepackage{amssymb}
\usepackage{amsfonts}

\begin{document}


\title[MAPPING PROPERTIES OF CERTAIN INTEGRAL OPERATORS]
      {MAPPING PROPERTIES OF CERTAIN INTEGRAL OPERATORS}

  \author{Victor Polunin}
\address{Chair of Applied Mathematics and Computer Modeling\\
    Belgorod State National Research University\\
         Pobedy street 85, Belgorod 308015, Russia}

        \email{polunin@bsu.edu.ru}

\author{Vladimir Vasilyev}
\address{Chair of Applied Mathematics and Computer Modeling\\
    Belgorod State National Research University\\
         Pobedy street 85, Belgorod 308015, Russia}

        \email{vbv57@inbox.ru}

 \author{Nelly Erygina}
\address{Chair of  Mathematics and Education in Natural Sciences\\
 Belgorod State National Research University\\
 Studencheskaya street 14, Belgorod 308007, Russia}

\email{erygina8@bsu.edu.ru}

\begin{abstract}
We study mapping properties of two-dimensional linear integral operators in some weighted spaces with special kernels. The considered spaces are certain variant of Sobolev--Slobodetskii spaces and their generalizations related to Banach spaces. Sufficient conditions for boundedness for such operators in these spaces are obtained.
\end{abstract}

\maketitle

\section{Introduction}

Some systems of integral equations arise under studying model pseudo-differential equations in plane corners \cite{V1,VV4}, for example when we consider simple boundary conditions like Dirichlet or Neumann conditions. According to this fact it seems interesting to study some mapping properties of such integral operators. These two-dimensional integral operator act in Fourier images of Sobolev--Slobodetski spaces which are weighted spaces with prescribed weight. These mapping properties play important role under studying discrete boundary value problems \cite{MVV,VVK} since we try to compare such integral operators with their truncated analogues, these arise under considering discrete boundary value problems. Let us note that integral operators and corresponding equations do not arise under studying model boundary value problems in a half-space, there are systems of linear algebraic equations as a result of reduction \cite{E}.

We will give here a simple example. Let $\Omega\subset\mathbb R^2$ be the first quadrant. We consider the following Dirichlet boundary value problem \cite{V1,VVK} in Sobolev--Slobodetskii space $H^s(\Omega)$
\begin{equation}\label{2}
\left\{
\begin{array}{rcl}
(Au)(x)=0,~~~x\in\Omega,\\
{u}_{|_{x_1=0}}=f(x_2),~~~~~~~~~~~~~{u}_{|_{x_2=0}}=g(x_1)
\end{array}
\right.
\end{equation}
where $A$ is a pseudo-differential operator with the symbol $A(\xi)$ satisfying the condition
\[
c_1(1+|\xi|)^{\alpha}\leq|A(\xi)|\leq c_2(1+|\xi|)^{\alpha}
\]
and admitting the wave factorization \cite{V1} with respect to $\Omega$
\[
A(\xi)=A_{\neq}(\xi)A_=(\xi)
\]
 with the index $\ae$ such that $\ae-s=1+\delta, |\delta|<1/2$ then the Dirichlet problem \eqref{2} can be reduced to the system of linear integral equations
\[
\left\{
\begin{array}{rcl}
&\int\limits_{-\infty}^{\infty}K_1(\xi)C(\xi_1)d\xi_1+D(\xi_2)=F(\xi_2)\\
&C(\xi_1)+\int\limits_{-\infty}^{\infty}K_2(\xi)D(\xi_2)d\xi_2=G(\xi_1)
\end{array}
\right.
\]
with respect to two unknown functions $C,D$. These functions should belong to Fourier image of certain $H^s$-space with  certain $s<0$.

\section{Integral operators}

We study the following integral operator
\[
(Kf)(x)=\int\limits_{-\infty}^{+\infty}K(x,y)f(y)dy
\]
in functional spaces with the following norm
\[
||f||_s=\left(\int\limits_{\mathbb R}|f(x)|^2(1+|x|)^{2s}dx\right)^{1/2}.
\]

We assume that the kernel $K(x,y)$ satisfies the following condition
\begin{equation}\label{1}
|K(x,y)|\sim(1+|x|+|y|)^{-\ae},~~~\ae\in\mathbb R.
\end{equation}

{\bf Theorem 1.} {\it
The operator $K$ is a linear bounded operator,
\[
K: {H}^{s_1}(\mathbb R)\rightarrow {H}^{s_2}(\mathbb R),
\]
where $s_1<0, \ae>\max\{1/2-s_1, 1+s_2-s_1\}$.
}

\begin{proof}
Let us verify.
\[
||Kf||^2_{s_2}=\int\limits_{-\infty}^{+\infty}(1+|x|)^{2s_2}|(Kf)(x)|^2dx=
\]
\[
=\int\limits_{-\infty}^{+\infty}(1+|x|)^{2s_2}\left|\int\limits_{-\infty}^{+\infty}K(x,y)f(y)dy\right|^2dx\leq
\]
\[
\leq\int\limits_{-\infty}^{+\infty}(1+|x|)^{2s_2}\left(\int\limits_{-\infty}^{+\infty}|K(x,y)|
|f(y)|dy\right)^2dx\leq
\]
\[
\leq~const\int\limits_{-\infty}^{+\infty}(1+|x|)^{2s_2}
\left(\int\limits_{-\infty}^{+\infty}(1+|x|+|y|)^{-\ae}|
|f(y)|dy\right)^2dx.
\]
In the inner integral, we apply the Cauchy--Schwartz inequality introducing the factors  $(1+|y|)^{-s_1}$ and $(1+|y|)^{s_1}$ form the first and the second term and taking into account that $(1+|y|)^{-s_1}\leq(1+|x|+|y|)^{-s_1}$. Then we obtain
\[
||Kf||^2_{s_2}\leq~cons||f||^2_{s_1}\int\limits_{-\infty}^{+\infty}(1+|x|)^{2s_2}
\left(\int\limits_{-\infty}^{+\infty}(1+|x|+|y|)^{-2(s_1+\ae)}
dy\right)dx\leq
\]
\[
\leq~const||f||^2_{s_1}\int\limits_{0}^{+\infty}(1+|x|)^{2s_2}
(1+|x|)^{-2(s_1+\ae-1/2)}
dx
\]
since $-2(s_1+\ae)<-1$. Thus, we conclude
\[
||Kf||^2_{s_2}\leq~cons||f||^2_{s_1}\int\limits_{0}^{+\infty}(1+|x|)^{2(s_2-s_1-\ae+1/2)}dx\leq~cons||f||^2_{s_1}
\]
since $2(s_2-s_1-\ae+1/2)<-1$.
\end{proof}

{\bf Corollary 1.} {\it
The operator $K$ is a linear bounded operator,
\[
K: {H}^{s}(\mathbb R)\rightarrow {H}^{s}(\mathbb R),
\]
where $s<0, \ae>\max\{1/2-s, 1\}$.
}

\section{Some generalizations: $H^{s,p}$-spaces}

\subsection{First variant}

This section is devoted to a generalization of the above result to Banach spaces $H^{s,p}(\mathbb R), 1<p<\infty$. Let us write that for $f\in H^{s,p}(\mathbb R)$
\[
||f||_{s,p}=\left(\int\limits_{\mathbb R}|f(x)|^p(1+|x|)^{ps}dx\right)^{1/p}.
\]
Obviously, for $p=2$ we have $H^{s,2}(\mathbb R)\equiv H^{s}(\mathbb R)$. Let us remind also that the pair $(p,q)$ is so that $1/p+1/q=1$.

{\bf Theorem 2.} {\it
Let the kernel $K(x,y)$ satisfies the condition \eqref{1}. Then the operator $K$ is a linear bounded operator
\[
K: {H}^{s_1,p_1}(\mathbb R)\rightarrow {H}^{s_2,p_2}(\mathbb R),
\]
where $s_1,s_2, \ae, p_1,p_2$ such that $s_1<0, \ae>\max\{1/{q_1}-s_1, 1/{p_2}+1/{q_1}+s_2-s_1\}$.
}

\begin{proof}
We estimate
\[
||Kf||^{p_2}_{s_2,p_2}=\int\limits_{-\infty}^{+\infty}(1+|x|)^{p_2s_2}|(Kf)(x)|^{p_2}dx=
\]
\[
=\int\limits_{-\infty}^{+\infty}(1+|x|)^{p_2s_2}\left|\int\limits_{-\infty}^{+\infty}K(x,y)f(y)dy\right|^{p_2}dx\leq
\]
\[
\leq\int\limits_{-\infty}^{+\infty}(1+|x|)^{p_2s_2}\left(\int\limits_{-\infty}^{+\infty}|K(x,y)|
|f(y)|dy\right)^{p_2}dx\leq
\]
\[
\leq~const\int\limits_{0}^{+\infty}(1+|x|)^{p_2s_2}
\left(\int\limits_{-\infty}^{+\infty}(1+|x|+|y|)^{-\ae}
|f(y)|dy\right)^{p_2}dx.
\]

Let us consider the inner integral
\[
\int\limits_{-\infty}^{+\infty}(1+|x|+|y|)^{-\ae}
|f(y)|dy=
\]
\[
=\int\limits_{-\infty}^{+\infty}(1+|y|)^{-s_1}(1+|x|+|y|)^{-\ae}
|f(y)|(1+|y|)^{s_1}dy
\]

Applying the H\"older inequality and taking into account the inequality $(1+|y|)^{-s_1}\leq(1+|x|+|y|)^{-s_1}$ we have
\[
\int\limits_{-\infty}^{+\infty}(1+|x|+|y|)^{-\ae}|
|f(y)|dy\leq||f||_{s_1,p_1}\left(\int\limits_{0}^{+\infty}(1+|x|+|y|)^{-q_1(s_1+\ae)}dy\right)^{\frac{1}{q_1}}\leq
\]
\[
\leq~const||f||_{s_1,p_1}(1+|x|)^{\frac{-q_1(s_1+\ae)+1}{q_1}},
\]
since $-q_1(s_1+\ae)<-1\Longleftrightarrow\ae>1/{q_1}-s_1$.

If so then we have
\[
||Kf||^{p_2}_{s_2,p_2}\leq~cons||f||^{p_2}_{s_1,p_1}\int\limits_{0}^{+\infty}(1+|x|)^{p_2(1/{q_1}-s_1-\ae+s_2)}
dx
\leq~const||f||^{p_2}_{s_1,p_1,}
\]
since $p_2(1/{q_1}-s_1-\ae+s_2)<-1\Longleftrightarrow\ae>1/{p_2}+1/{q_1}+s_2-s_1$.
\end{proof}

{\bf Corollary 2.} {\it
Let the kernel $K(x,y)$ satisfies the condition \eqref{1}. Then the operator $K$ is a linear bounded operator
\[
K: {H}^{s,p}(\mathbb R)\rightarrow {H}^{s,p}(\mathbb R),
\]
where $s, \ae, p$ such that $s<0, \ae>\max\{1/{q}-s, 1\}$.
}

\subsection{Second variant}

This section is devoted to another generalization of the above result to Banach spaces $H^{s,p}(\mathbb R), 1<p<\infty$. We say that for $f\in H^{p,s}(\mathbb R)$
\[
||f||_{p,s}=\left(\int\limits_{\mathbb R}|f(x)|^p(1+|x|)^{2s}dx\right)^{1/p}.
\]
Obviously, for $p=2$ we have $H^{s,2}(\mathbb R)\equiv H^{s}(\mathbb R)$. Let us remind also that the pair $(p,q)$ is so that $1/p+1/q=1$.

{\bf Theorem 3.} {\it
Let the kernel $K(x,y)$ satisfies the condition \eqref{1}. Then the operator $K$ is a linear bounded operator
\[
K: {H}^{p_1,s_1}(\mathbb R)\rightarrow {H}^{p_2,s_2}(\mathbb R),
\]
where $s_1,s_2, \ae, p_1,p_2$ such that $s_1<0, \ae>\max\{1/{q_1}-2s_1/{p_2}, 1/{p_2}+1/{q_1}+2s_2/{p_2}-2s_1/{p_1}\}$.
}

\begin{proof}
Let us estimate the norm
\[
||Kf||^{p_2}_{p_2,s_2}=\int\limits_{-\infty}^{+\infty}(1+|x|)^{2s_2}|(Kf)(x)|^{p_2}dx=
\]
\[
=\int\limits_{-\infty}^{+\infty}(1+|x|)^{2s_2}\left|\int\limits_{-\infty}^{+\infty}K(x,y)f(y)dy\right|^{p_2}dx\leq
\]
\[
\leq\int\limits_{-\infty}^{+\infty}(1+|x|)^{2s_2}\left(\int\limits_{-\infty}^{+\infty}|K(x,y)|
|f(y)|dy\right)^{p_2}dx\leq
\]
\[
\leq~const\int\limits_{0}^{+\infty}(1+|x|)^{2s_2}
\left(\int\limits_{-\infty}^{+\infty}(1+|x|+|y|)^{-\ae}
|f(y)|dy\right)^{p_2}dx.
\]

Let us consider the inner integral
\[
\int\limits_{-\infty}^{+\infty}(1+|x|+|y|)^{-\ae}
|f(y)|dy=
\]
\[
=\int\limits_{-\infty}^{+\infty}(1+|y|)^{-2s_1/{p_1}}(1+|x|+|y|)^{-\ae}
|f(y)|(1+|y|)^{2s_1/{p_1}}dy
\]

Applying the H\"older inequality and taking into account the inequality $(1+|y|)^{-2s_1/{p-1}}\leq(1+|x|+|y|)^{-2s_1/{p-1}}$ we have
\[
\int\limits_{-\infty}^{+\infty}(1+|x|+|y|)^{-\ae}|
|f(y)|dy\leq||f||_{p_1,s_1}\left(\int\limits_{0}^{+\infty}(1+|x|+|y|)^{-q_1(2s_1/{p_1}+\ae)}dy\right)^{\frac{1}{q_1}}\leq
\]
\[
\leq~const||f||_{p_1,s_1}(1+|x|)^{\frac{-q_1(2s_1/{p_1}+\ae)+1}{q_1}},
\]
since $-q_1(2s_1/{p_1}+\ae)<-1\Longleftrightarrow\ae>1/{q_1}-2s_1/{p_2}$.

If so then we have
\[
||Kf||^{p_2}_{p_2,s_2}\leq~cons||f||^{p_2}_{p_1,s_1}\int\limits_{0}^{+\infty}(1+|x|)^{p_2(1/{q_1}-2s_1/{p_1}-\ae)+2s_2}
dx
\leq~const||f||^{p_2}_{p_1,s_1,}
\]
since $p_2(1/{q_1}-2s_1/{p_1}-\ae)+2s_2<-1\Longleftrightarrow\ae>1/{p_2}+1/{q_1}+2s_2/{p_2}-2s_1/{p_1}$.
\end{proof}

{\bf Corollary 3.} {\it
Let the kernel $K(x,y)$ satisfies the condition \eqref{1}. Then the operator $K$ is a linear bounded operator
\[
K: {H}^{p,s}(\mathbb R)\rightarrow {H}^{p,s}(\mathbb R),
\]
where $s, \ae, p$ such that $s<0, \ae>\max\{1/{q}-2s/{p}, 1\}$.
}

\section*{Conclusion}

The considered spaces $H^{s,p}(\mathbb R)$ and $H^{p,s}(\mathbb R)$ like almost  the analogous spaces considered in \cite{T,MW}. But we work in Fourier images, and the introduced spaces are usual weighted spaces with power weights.



\begin{thebibliography}{99}

\bibitem{V1}  V.~B.~Vasil'ev, Wave
 Factorization of Elliptic Symbols: Theory and
Applications. Introduction to the Theory of Boundary Value Problems in Non-Smooth Domains.
Dordrecht--Boston--London: Kluwer Academic Publishers, 2000.

\bibitem{VV4} A.~V.~Vasilyev, V.~B.~Vasilyev, Difference equations and boundary value problems.  ~In: Pinelas, S., D\u{o}sl\'a, Z., D\u{o}sl\'y, O., Kloeden, P. (eds.) Differential and Difference Equations and Applications, Springer Proc. Math. \& Stat. 2016, \textbf{164}  421--432.

    \bibitem{E} G.~Eskin,
Boundary Value Problems for Elliptic Pseudodifferential Equations.
Providence, RI: AMS, 1981.

 \bibitem{MVV} A.~A.~Mashinets, A.~V.~Vasilyev, V.~B.~Vasilyev, On discrete Neumann problem in a quadrant. Lobachevskii J. Math. 2023, \textbf{44}~(3) 1011--1021.


 \bibitem{VVK} A.~V.~Vasilyev, V.~B.~Vasilyev, A.`A.~Khodyrevz, The discrete Dirichlet problem. Solvability and approximation properties. J. Math. Sci. 2023, \textbf{270}~(5) 654--664.

 \bibitem{T} H.~Triebel, Theory of Function Spaces.
Basel--Boston--Stuttgart: Birkh\"auser Verlag, 1983.

 \bibitem{MW} M.W.~Wong, An Introduction to Pseudo-Differential Operators.
 Singapore--New Jersey--London--Hong Kong: World Scientific, 1991.

 \end{thebibliography}
\end{document}